\documentclass[10pt]{amsart}
\usepackage{amssymb, latexsym}
\usepackage{enumerate}

\theoremstyle{plain}
\newtheorem{thm}{Theorem}

\newtheorem{prop}[thm]{Proposition}
\newtheorem{lem}[thm]{Lemma}

\numberwithin{equation}{section} \numberwithin{thm}{section}

\begin{document}


\title[Bootstrapped Morawetz Estimates and resonant decompositions for
NLS]{Bootstrapped Morawetz estimates and resonant decomposition for
low regularity global solutions of cubic NLS on $\mathbb{R}^2$}
\author{J. Colliander}
\thanks{J.C. was supported in part by NSERC grant RGP250233-07.}
\address{Department of Mathematics, University of Toronto, Toronto, ON, Canada M5S 2E4}
\email{\tt colliand@math.toronto.edu}

\author{T. Roy}
\address{Department of Mathematics, University of California, Los Angeles} \email{\tt triroy@math.ucla.edu}

\date{9 November 2008}




\vspace{-0.3in}

\begin{abstract}
We prove global well-posedness for the $L^{2}$-critical cubic
defocusing nonlinear Schr\"odinger equation on $\mathbb{R}^{2}$ with
data $u_{0} \in H^{s}(\mathbb{R}^{2})$ for $ s > \frac{1}{3}$. The
proof combines  \textit{a priori} Morawetz estimates obtained in
\cite{cgt} and the improved almost conservation law obtained in
\cite{cksttI2}. There are two technical difficulties. The first one
is to estimate the variation of the improved almost conservation law
on intervals given in terms of Strichartz spaces rather than in terms of $X^{s,b}$ spaces. The
second one is to control the error of the \textit{a priori} Morawetz
estimates on an arbitrary large time interval, which is performed by
a bootstrap via a double layer in time decomposition.
\end{abstract}
\maketitle

\section{Introduction}

We shall consider the $L^{2}$-critical Schr\"odinger equation on
$\mathbb{R}^{2}$

\begin{equation}
\begin{array}{ll}
iu_{t}  + \Delta u  &  =   |u|^{2}u
\end{array}
\label{Eqn:Schroddefcubic}
\end{equation}
with data $u(0)=u_{0} \in H^{s}(\mathbb{R}^{2})$, $s \geq 0$. Here
$H^{s}(\mathbb{R}^{2})$ denotes the Sobolev space endowed with the
norm

\begin{equation}
\begin{array}{ll}
\| f \|_{H^{s}(\mathbb{R}^{2})} & := \|<\xi>^{s} \widehat{f}(\xi)
\|_{L^{2}(\mathbb{R}^{2})}
\end{array}
\end{equation}
with $\widehat{f}$ denoting the Fourier transform

\begin{equation}
\begin{array}{ll}
\widehat{f}(\xi) & := \int_{\mathbb{R}^{2}} f(x) e^{-ix \cdot \xi}
\, dx
\end{array}
\end{equation}
and $\langle \xi \rangle:= \left( 1 + |\xi|^{2}  \right)^{\frac{1}{2}}$. This
problem is known to be locally well-posed \cite{cazwei} for any $s
\geq 0$. If $s>0$ then local well-posedness means that for any data
$u_{0} \in H^{s} (\mathbb{R}^{2})$, there exists a time of local
existence $T_{l}= T_{l} (\| u_{0} \|_{H^{s}(\mathbb{R}^{2})})$
depending only on the norm of the initial data and a unique solution
$u$ lying in a Banach space $X \subset \mathcal{C} \left( [0,T_{l}],
H^{s}(\mathbb{R}^{2}) \right)$ such that $u(t)$ satisfies for $t \in
[0,T_{l}]$ the Duhamel formula

\begin{equation}
\begin{array}{ll}
u(t) : = e^{it \triangle} u_{0} - i \int_{0}^{t} e^{i(t-t^{'})
\triangle} \left[ |u|^{2} u (t^{'}) \right] \, dt^{'}
\end{array}
\end{equation}
and the solution depends continuously on the norm of the initial
data. Local-in-time $H^{s}$-solutions to (\ref{Eqn:Schroddefcubic})
satisfy the \textit{ mass conservation law}

\begin{equation}
\begin{array}{ll}
\| u(t) \|_{L^{2}(\mathbb{R}^{2})} & = \| u_{0} \|_{\mathbb{R}^{2}}
\end{array}
\end{equation}
and local-in-time $H^{1}$-solutions to (\ref{Eqn:Schroddefcubic})
satisfy the \textit{energy conservation law}

\begin{equation}
\begin{array}{ll}
E \left(u (t) \right) & : = \frac{1}{2} \int_{\mathbb{R}^{2}}
|\nabla u(t,x)|^{2} \, dx + \frac{1}{4} \int_{\mathbb{R}^{2}}
|u(t,x)|^{4} \, dx.
\end{array}
\label{Eqn:NrjConservLaw}
\end{equation}

In this paper we are interested in proving that $H^{s}$-solutions to
(\ref{Eqn:Schroddefcubic}) with $s \geq 0$ exist for all time $T
\geq 0$. If $s >0$ then in view of the local well-posedness theory
it suffices to find an \textit{a priori} bound of the form

\begin{equation}
\begin{array}{ll}
\| u(T) \|_{H^{s}(\mathbb{R}^{2})} & \leq Q \left( \| u_{0}
\|_{H^{s}(\mathbb{R}^{2})}, T \right)
\end{array}
\label{Eqn:Bound}
\end{equation}
with $Q$ a function depending only on the norm of the initial data
and time $T$. If $s=1$ then the energy conservation law immediately
yields the bound (\ref{Eqn:Bound}). No blowup solutions are known
for (\ref{Eqn:Schroddefcubic}). It is conjectured that
(\ref{Eqn:Schroddefcubic}) is globally well-posed in
$H^{s}(\mathbb{R}^{2})$ for $1> s \geq 0$.
The first breakthrough to establish global well-posedness below the energy
threshold, by using what is now referred to as the \textit{Fourier
truncation method}, appears in \cite{bourg}. He showed global well-posedness for data in
$H^{s}(\mathbb{R}^{2})$ with $s
> \frac{3}{5}$. A sequence of works (\cite{cksttI1,cksttI2,fgr,cgt})
has lowered the regularity requirements for global well-posedness
for (\ref{Eqn:Schroddefcubic}) down to $s>\frac{2}{5}$. Recently, the conjecture was proved in \cite{ktv},
in the case of spherically symmetric initial data. The main result
of this paper is the following improvement:

\begin{thm}
The $L^{2}$-critical Schr\"odinger equation on $\mathbb{R}^{2}$ is
globally well-posed in $H^{s}(\mathbb{R}^{2})$, $1 > s > \frac{1}{3}
$. Moreover there exists a constant $C$ depending only on $\| u_{0}
\|_{H^{s}(\mathbb{R}^{2})}$ such that

\begin{equation}
\begin{array}{ll}
\| u(T) \|^{2}_{H^{s}(\mathbb{R}^{2})} & \leq C( \| u_{0}
\|_{H^{s}(\mathbb{R}^{2})}) T ^{\frac{1-s}{3s-1}+}
\end{array}
\label{Eqn:BoundHsuT}
\end{equation}
for all times $T$.
\label{Thm:Gwp13}
\end{thm}

Before sketching the main ideas underpinning this theorem, we set up
some notation.

Given $A,B$ two nonnegative numbers, $A \lesssim B$ means that there
exists a universal nonnegative constant $K$ such that $A \leq K B$.
We say that $K_{0}$ is the constant determined by the relation $A
\lesssim B$ if $K_{0}$ is the smallest $K$ such that $A \leq K B$ is
true. We write $A \sim B$ when $A \lesssim B$ and $B \lesssim A$. $A
<< B$ denotes $ A \leq K B$ for some universal constant $K <
\frac{1}{100}$ . We also use the notations $A+ = A + \epsilon$,
$A++=A + 2 \epsilon$, $A- =A - \epsilon$ and $A-- = A - 2 \epsilon$,
etc. for some universal constant $0 < \epsilon \ll 1$. We shall
abuse the notation and write $+$, $-$ for $0+$, $0-$ respectively.

Let $ \lambda \in \mathbb{R}$ and let $u^{\lambda}$ denote the
following function

\begin{equation}
\begin{array}{ll}
u^{\lambda}(t,x) & := \frac{1}{\lambda} u \left(
\frac{t}{\lambda^{2}}, \frac{x}{\lambda} \right).
\end{array}
\label{Eqn:Defulambda}
\end{equation}
We recall that if $u$ satisfies (\ref{Eqn:Schroddefcubic}) with data
$u_{0}$ then $u^{\lambda}$ also satisfies (\ref{Eqn:Schroddefcubic})
but with data $ \frac{1}{\lambda} u_{0} \left( \frac{x}{\lambda}
\right)$.

If $J:=[a,b]$ is an interval then $|J|$ is its size. A partition
$\mathcal{P}_{\mu}(J)=(J_{i})_{i \in \{1, ..., l \}}$ of a finite
interval $J$ is of size $\mu$, $\mu>0$ if three conditions are
satisfied

\begin{enumerate}

\item $ \bigcup_{i \in [1,..l]} J_{i} = J $

\item $J_{i} \cap J_{j} = \emptyset, i \neq j $

\item $|J_{i}|= \mu, i \in \{1, ..., l-1 \}$.

\end{enumerate}
If $u$ if a solution of (\ref{Eqn:Schroddefcubic}) on $J$ then we
can write $u$ as the sum of its linear part and its nonlinear part;
more precisely

\begin{equation}
\begin{array}{ll}
u(t) & = u_{J}^{l}(t) + u_{J}^{nl}(t)
\end{array}
\end{equation}
with

\begin{equation}
\begin{array}{ll}
u_{J}^{l}(t) : = e^{i (t-a) \triangle} u(a)
\end{array}
\end{equation}
and

\begin{equation}
\begin{array}{ll}
u_{J}^{nl}(t) := - i \int_{a}^{t} e^{i(t-t^{'}) \triangle} \left[
|u|^{2} u (t^{'}) \right] \, dt^{'}.
\end{array}
\end{equation}
Let $dt$ denote the standard Lebesgue measure and let $d \mu_{l}$
\footnote{$\mu_{l}$ since this measure will be applied to be the
linear part $u_{J}^{l}$ of $u$} be the following measure

\begin{equation}
\begin{array}{ll}
d \mu_{l} & : = \delta(t-a) dt.
\end{array}
\end{equation}
If $ (p,q) \in [1, \infty]$ then we define the spaces $L^{p}(J)$ and
$ L^{p}(J, d \mu^{l})$

\begin{equation}
\begin{array}{ll}
L^{p} (J) & := \left\{ f:\mathbb{R} \rightarrow \mathbb{C}, \, \| f
\|^{p}_{L^{p}}
:= \int_{J} |f|^{p} \, dt < \infty \right\} \\
L^{p}(J, d \mu^{l}) & := \left\{ f:\mathbb{R} \rightarrow
\mathbb{C}, \, \| f \|^{p}_{L^{p}( d \mu^{l})} := \int_{J} |f|^{p}
\, d \mu_{l} < \infty  \right\}
\end{array}
\end{equation}
and the mixed spaces

\begin{equation}
\begin{array}{ll}
L_{t}^{p}(J) L_{x}^{q} & := \left\{ f:\mathbb{R}^{2+1} \rightarrow
\mathbb{C}, \| f \|^{p}_{L_{t}^{p}(J) L_{x}^{q}} := \int_{J} \left(
\int_{\mathbb{R}^{2}} |f(t,x)|^{q} \, dx \right)^{\frac{q}{p}}
dt < \infty \right\} \\
L_{t}^{p}(J, d \mu_{l}) L_{x}^{q} & := \left\{ f:\mathbb{R}^{2+1}
\rightarrow \mathbb{C}, \, \| f \|^{p}_{L_{t}^{p}(J, d \mu_{l})
L_{x}^{q}} := \int_{J}  \left( \int_{\mathbb{R}^{2}} |f(t,x)|^{q} \,
dx \right)^{\frac{q}{p}} d \mu_{l} < \infty \right\}
\end{array}
\end{equation}
Let $\tilde{f}$ be the spacetime Fourier transform of a function $f$

\begin{equation}
\begin{array}{ll}
\tilde{f}(t,x) & : = \int_{\mathbb{R}^{2+1}} f(t,x) e^{-i(t \tau + x
\xi)} \, dt \, dx
\end{array}
\end{equation}
If $p$ is an integer larger or equal to one, $\sigma:
\mathbb{R}^{2p} \rightarrow \mathbb{C} $ is a smooth symbol and
$u_{1}$,...,$u_{2p}$ are Schwartz functions then we define the $2p$-linear functionals

\begin{equation}
\begin{array}{ll}
\Lambda_{2p}(\sigma;u_{1}(t),...,u_{2p}(t)):=
\int_{\xi_{1}+...+\xi_{2p}=0} \sigma(\xi_{1},...,\xi_{2p}) \prod_{ j
\, odd } \widehat{u_{j}}(t,\xi_{j}) \prod_{j \, even}
\widehat{\overline{u}_{j}}(t,\xi_{j})
\end{array}
\end{equation}
and

\begin{equation}
\begin{array}{ll}
\Lambda_{2p,J}(\sigma;u_{1},...,u_{2p}):=  \int_{J}
\int_{\xi_{1}+...+\xi_{2p}=0} \sigma(\xi_{1},...,\xi_{2p}) \prod_{ j
\, odd } \widehat{u_{j}}(t,\xi_{j}) \prod_{j \, even}
\widehat{\overline{u}_{j}}(t,\xi_{j}).
\end{array}
\end{equation}
If $u_{1}=...=u_{2p}=u$ then we abbreviate $\Lambda_{2p}(\sigma;u)
:= \Lambda_{2p}(\sigma;u_{1},..,u_{2p})$ and
$\Lambda_{2p,J}(\sigma;u) :=
\Lambda_{2p,J}(\sigma;u_{1},..,u_{2p})$. Let $\Omega^{1 \rightarrow
2p}_{k}$ denote the set of unordered subsets of size $k$ from the
set $\{0,...,2p \}$. If $A \in \Omega_{k}^{1 \rightarrow 2p}$ then
we write $\Lambda_{2p,J,A}(\sigma;u)$ for
$\Lambda_{2p,J}(\sigma;v_{1},...,v_{2p})$ with $v_{i}=u_{J}^{l}$ if
$i \in A$ and $v_{i}=u_{J}^{nl}$ if $i \notin A$. Let

\begin{equation}
\begin{array}{ll}
L &: = \bigcap_{i=1}^{4} \left\{ (\xi_{j})_{j \in [1,..,4 ]}, \,
|\xi_{j}| \leq \frac{N}{100}   \right\}
\end{array}
\end{equation}
and

\begin{equation}
\begin{array}{ll}
\Gamma & : = \left\{ (\xi_{j})_{j \in [1,..,4]}, \, \left|
\cos{(\xi_{12},\xi_{14})} \right| \geq \theta \right\}
\end{array}
\end{equation}
where $0 < \theta /ll 1$ is a parameter to be determined. Here we use
the convention $\xi_{ab}:=\xi_{a} + \xi_{b}$, $\xi_{abc}:= \xi_{a} +
\xi_{b} + \xi_{c}$, etc.

We constantly use the $I$- method \cite{cksttI1} throughout this
paper in order to find a pointwise-in-time upper bound of the
$H^{s}$-norm of the solution to (\ref{Eqn:Schroddefcubic}) with data
$u(0)=u_{0} \in H^{s}(\mathbb{R}^{2})$. We recall it now. Let $I$ be
the following multiplier

\begin{equation}
\begin{array}{ll}
\widehat{If}(\xi) & := m(\xi) \hat{f}(\xi)
\end{array}
\end{equation}
where $m(\xi): =  \eta \left( \frac{\xi}{N} \right)$, $\eta$ is a
smooth, radial, nonincreasing in $|\xi|$ such that

\begin{equation}
\begin{array}{ll}
\eta (\xi) & :=  \left\{
\begin{array}{l}
1, \, |\xi| \leq 1 \\
\frac{1}{|\xi|^{1-s}}, \, |\xi| \geq 2
\end{array}
\right.
\end{array}
\end{equation}
By plugging the multiplier $I$ into the energy conservation law
(\ref{Eqn:NrjConservLaw}) we define the so-called modified energy

\begin{equation}
\begin{array}{ll}
E \left( I u(t) \right) & : = \frac{1}{2} \int_{\mathbb{R}^{2}}
|\nabla I u(t,x)|^{2}  dx + \frac{1}{4} \int_{\mathbb{R}^{2}} | Iu
(t,x)|^{4} \, dx.
\end{array}
\end{equation}
The following proposition \cite{cksttI1} shows that it suffices to
estimate the modified energy at time $T$ in order to find an upper
bound of the pointwise-in-time $H^{s}$-norm of the solution $u$ to
(\ref{Eqn:Schroddefcubic}) with data $u(0)=u_{0} \in
H^{s}(\mathbb{R}^{2})$; more precisely

\begin{prop}[$H^{s}(\mathbb{R}^{2})$ norm and modified energy are comparable \cite{cksttI1}]
For all time $T \geq 0$
\begin{equation}
\begin{array}{ll}
\| u(T) \|^{2}_{H^{s}(\mathbb{R}^{2})} & \lesssim E(Iu(T)) + \|
u_{0} \|^{2}_{H^{s}(\mathbb{R}^{2})}.
\end{array}
\label{Eqn:HsEIu}
\end{equation}
\label{prop:HsnormNrj}
\end{prop}
Since the symbol of $I$ approaches one as $N$ goes to infinity we
expect the variation of the modified energy to be slower and slower
as $N$ increases. Therefore we estimate the modified energy by
using the fundamental theorem of calculus and we use Proposition
\ref{prop:HsnormNrj} to control  $\| u(T) \|_{H^{s}(\mathbb{R}^{2})}
$. \vspace{5 mm}

The paper is organized as follows:

In Section
\ref{sec:NewAlmConsLaw} we recall the main ideas of \cite{cksttI2}.
In particular, we explain their construction of a new almost
conservation law $\tilde{E}(u(t))$  which is close to the modified
energy $E(Iu(t))$ at each time $t$ and how they estimated the
variation of $\tilde{E}(u(t))$ on an interval of size, roughly
speaking, equal to one. In Section \ref{sec:Morawetz} we recall the
main results of \cite{cgt} and, in particular, the Morawetz-type
estimates. We would like to combine the ideas from \cite{cksttI2}
with those from \cite{cgt}. However there are two non-trivial
difficulties that appear.

The $I$-method is based upon an
estimate of the variation of an almost conservation law on a small
interval where we have a control of a large number of norms. Then
the variation of the almost conservation law on an arbitrary large
time interval $[0,T]$ is estimated by iteration on each subinterval
of a partition of $[0,T]$ where this local control holds. This total
variation must be controlled at the end of the process. Therefore,
if we can establish a local control on a subinterval as large as
possible then the number of iterations is reduced and we have a
better control of the total variation, which implies global
well-posedness for rougher data. Unfortunately we cannot use the
result established in \cite{cksttI2} (see Proposition
\ref{prop:PropAlmCons}) to estimate the variation of
$\tilde{E}(u(t))$ since the local control of the solution in
$X^{s,b}$ spaces is only true for short time intervals (see
Proposition \ref{Prop:ModLocWp}). This is due to the nature of these
spaces: they describe very well the solution locally in time but not
on long time intervals. Proposition \ref{prop:LocalWellPosedness}
shows that we have a local control on intervals $J$ where the
$L_{t}^{4} L_{x}^{4}$ norm of $Iu$ is small. The first idea would be
to divide $[0,T]$ into subintervals $J$ where the $L_{t}^{4}
L_{x}^{4}$ norm of $Iu$ is small. Indeed their size is expected to
be, roughly speaking, larger than one, because  the Morawetz-type
estimates provide good control of the $L_{t}^{4} L_{x}^{4}$ norm of
$Iu$ on $[0,T]$. In Section \ref{sec:AlmCon} we estimate the
variation of $\tilde{E}$ on $J$. The proof has similarities with
that of (\ref{Eqn:AlmConsXsb}) but there are differences in the
method. We write the variation of $\tilde{E}$ on $J$ in the
spacetime Fourier domain, we decompose $u$ into the sum of its
linear part $u^{l}_{J}$ and its nonlinear part $u^{nl}_{J}$, and
after some measure rearrangements performed via the use of Fubini's
theorem, we use some refined bilinear estimates
\cite{bourg,cksttI2}. These estimates are key estimates to get a
slow increase of $\tilde{E}$. At the end of the process we can bound
the variation of $\tilde{E}$ by some quantities that are estimated
by the local control theory in turn.

Unfortunately, if we use the Morawetz-type estimate on the whole
$[0,T]$ then an error term appears and, as time $T$ goes to
infinity, it grows at a faster rate than that generated by the
variation of the modified energy on the same interval. The control
of the error term is possible if and only if $s>\frac{2}{5}$ (see
\cite{cgt}). We would like to use the Morawetz-type estimates in a
better way. To this end we perform a double layer in time
decomposition. First we divide $[0,T]$ into subintervals $J$ of
size, roughly speaking, equal to $ N^{3-}$. This enables us to
control the error term of the Morawetz estimate on $J$ by its main
term. Then we decompose each $J$ into subintervals $J_{k}$ where the
$L_{t}^{4} L_{x}^{4}$ norm of $ I u $ is small. By applying the
local control theory and the almost conservation law in $L_{t}^{q}
L_{x}^{r}$ spaces (see Proposition \ref{prop:VarNewAlmStr}) we can
estimate the variation of $\tilde{E}$ on $J_{k}$ and then on $J$ by
iteration. The final step is to bootstrap the Morawetz estimates.
More precisely, we use for every $J$ the corresponding Morawetz
estimate and we iterate to estimate the variation of $\tilde{E}$ on
the whole interval $[0,T]$. At the end of the proof we can control
the modified energy on $[0,T]$, provided that $s>\frac{1}{3}$. The
whole process is explained in Section \ref{sec:PfGw13}.

\vspace{5mm}

$\textbf{Acknowledgements}:$ We would like to thank T.~Tao and N.~Tzirakis for
discussions related to this work.

\section{Summary of \cite{cksttI2} }
\label{sec:NewAlmConsLaw}

In this section we recall the main ideas and results of
\cite{cksttI2} since we will often refer to them throughout this
paper.

The variation of the modified energy $E(Iu(t))$ is not equal to
zero, because of the presence of the commutator $I u^{3} -
(Iu)^{3}$. The control of this variation is possible if the Sobolev
exponent $s$ is larger than a threshold $s_{0}$. This control
implies global existence for data in $H^{s}(\mathbb{R}^{2})$,
$s>s_{0}$. In \cite{cksttI2} the authors aimed at designing a new
almost conservation law $\tilde{E}(u(t))$ that would satisfy two
properties

\begin{enumerate}

\item Almost conservation law: $\tilde{E}(u(t))$ would have a slower
variation than $E(Iu(t))$

\item Proximity to $E(Iu(t))$ at each time $t$: this property would
allow to control $E(Iu(t))$  \emph{via} $\tilde{E}(u(t))$

\end{enumerate}
To this end they searched for a candidate $\tilde{E}$ that would
have the following form

\begin{equation}
\begin{array}{ll}
\widetilde{E}(u(t)) & := \frac{1}{2} \Lambda_{2}(\sigma_{2};u(t)) +
\Lambda_{4}(\sigma_{4};u(t))
\end{array}
\label{Eqn:DefTildeNrj}
\end{equation}
with $\sigma_{2}$ denoting the following multiplier

\begin{equation}
\begin{array}{ll}
\sigma_{2} & := - \xi_{1} m(\xi_{1}). \xi_{2} m(\xi_{2})
\end{array}
\end{equation}
and $\sigma_{4}$ to be determined. Notice that
$\Lambda_{2}(\sigma_{2})(u(t))$ is nothing else but the kinetic part
of the modified energy, i.e $\Lambda_{2}(\sigma_{2})(u(t))=
\frac{1}{2} \| I u(t) \|^{2}_{\dot{H}^{1}}$. The idea is to
substitute the potential term $ V(t):= \frac{1}{4}
\int_{\mathbb{R}^{2}} |Iu(t,x)|^{4} \, dx$ of the modified energy
$E(Iu)$ for a new quadrilinear term $ \Lambda_{4}(\sigma_{4};u(t))$
and to search for some cancellations in the computation of the
derivative $\partial_{t} \tilde{E}(u(t))$. If we compute the
derivative of $\Lambda_{2}(\sigma_{2};u(t))$ and
$\Lambda_{4}(\sigma_{4};u(t))$ then we find, by using
(\ref{Eqn:Schroddefcubic})

\begin{equation}
\begin{array}{ll}
\partial_{t} \left( \frac{1}{2} \Lambda_{2} (\sigma_{2};u(t)) \right) :& =  \Lambda_{4} (\mu; u (t))
\end{array}
\end{equation}
and

\begin{equation}
\begin{array}{ll}
\partial_{t} \Lambda_{4} ( \sigma_{4} \alpha_{4};u(t)) & : = \Lambda_{4} (\alpha_{4} \sigma_{4};u(t)) +
\Lambda_{6} ( \nu_{6};u(t))
\end{array}
\end{equation}
with

\begin{equation}
\begin{array}{ll}
\mu & := \frac{i}{4} \left( |\xi_{1}|^{2}m^{2}(\xi_{1}) -
|\xi_{2}|^{2}m^{2}(\xi_{2}) + |\xi_{3}|^{2} m^{2}(\xi_{3}) -
|\xi_{4}|^{2} m^{2}(\xi_{4} ) \right)
\end{array}
\label{Eqn:Defmu}
\end{equation}

\begin{equation}
\begin{array}{ll}
\alpha_{4} & := -i \left( |\xi_{1}|^{2} - |\xi_{2}|^{2} +
|\xi_{3}|^{2} - |\xi_{4}|^{2} \right)
\end{array}
\label{Eqn:Defalpha4}
\end{equation}
and

\begin{equation}
\begin{array}{ll}
\nu_{6} & := -i \sum_{k=1}^{4} (-1)^{k+1}
\sigma_{4}(\xi_{1},...,\xi_{k}+..+\xi_{k+2},\xi_{k+3},...,\xi_{6})
\end{array}
\end{equation}
The authors tried to cancel the quadrilinear terms resulting from
the derivative of $\widetilde{E}(u(t))$ by letting $\sigma_{4}:= -
\frac{\mu}{\alpha_{4}}$. The problem is that the singularity
$\alpha_{4}=0$ appears. Therefore they had to truncate optimally
$\sigma_{4}$ away from $\alpha_{4}=0$  so that the truncation does
not totally lose the effect of these cancellations. This requires a
detailed study of the singularity. Recall that the corrective term
$\Lambda_{4}(\sigma_{4};u(t))$ is a quadrilinear integral evaluated
on the convolution surface $\xi_{1}+...+\xi_{4}=0$. They observed
that $\alpha_{4}=2i \xi_{12} \xi_{14} \cos{(\xi_{12}, \xi_{14})}$ on
this surface and that if $|\xi_{i}| \ll N$, $i \in \{1,..,4 \}$ then
the singularity disappears and $\alpha_{4}=\frac{1}{4}$
\footnote{This phenomenon is expected. Indeed if all the frequencies
have amplitude smaller than $N$ then the modified energy is the
energy itself and the variation is equal to zero.}. Therefore they
truncated $ -\frac{\mu}{\alpha_{4}}$ in the following way

\begin{equation}
\begin{array}{ll}
\sigma_{4}(\xi_{1},...,\xi_{4}) & : = - \frac{\mu}{\alpha_{4}}
\chi_{L \cup \Gamma} (\xi_{1},...,\xi_{4})
\end{array}
\label{Eqn:Defsigma4}
\end{equation}
With this value for $\sigma_{4}$, $\widetilde{E}(u(t))$ is
well-defined by (\ref{Eqn:DefTildeNrj}). They showed that
$\tilde{E}(u(t))$ and $E(Iu(t))$ are closed to each other at each
time $t$; more precisely

\begin{prop}[Proximity to E(Iu(t)) at each time $t$  \cite{cksttI2}]
\begin{equation}
\begin{array}{ll}
\left| \widetilde{E}\left( u(t) \right)) - E \left( Iu(t) \right)
\right| \lesssim \frac{1}{\theta N^{2-}} \| I
u(t)\|^{4}_{H^{1}(\mathbb{R}^{2})}.
\end{array}
\label{Eqn:FixedEst}
\end{equation}
\label{prop:FixedEst}
\end{prop}
Then, by using a delicate mutilinear analysis, they proved the
following result

\begin{prop}[Almost Conservation Law in $X^{s,b}$ spaces \cite{cksttI2}]
\begin{equation}
\begin{array}{ll}
\left| \sup_{t \in J} \widetilde{E}(u(t)) - \widetilde{E}(u(a))
\right| & \lesssim \left(
\frac{\theta^{\frac{1}{2}}}{N^{\frac{3}{2}-}} + \frac{1}{N^{2-}} +
\frac{1}{\theta N^{3-}} \right) \| I u \|^{4}_{X^{1,\frac{1}{2}+}}
\end{array}
\label{Eqn:AlmConsXsb}
\end{equation}
\label{prop:PropAlmCons}
\end{prop}
The definition of the $X^{1,\frac{1}{2}+}$ spaces can be found in
\cite{bourg2} for example. The proof of Proposition
\ref{prop:PropAlmCons} extensively relies upon two refined bilinear
estimates

\begin{prop}[Bilinear estimates \cite{bourg}, \cite{cksttI2}] Let $f$, $g$ be two Schwartz
functions. Let $N_{1}$, $N_{2}$ be two dyadic numbers such that
$N_{1}<N_{2}$. Let $\theta$ be a parameter such that $0< \theta \ll
1$. If

\begin{equation}
\begin{array}{ll}
B_{\epsilon}(\tau,\xi) :=  \int_{\xi_{1}+\xi_{2}=\xi}
\chi_{_{|\xi_{1}| \sim N_{1}}} \chi_{_{|\xi_{2}| \sim N_{2}}}
\chi_{_{[-\epsilon,\epsilon]}} (\tau - |\xi_{1}|^{2} -
|\xi_{2}|^{2}) \widehat{f}(\xi_{1}) \widehat{g}(\xi_{2}) d \xi_{1}
\end{array}
\label{Eqn:Def1Bil}
\end{equation}
and

\begin{equation}
\begin{array}{ll}
B_{\epsilon,\theta}(\tau,\xi) :=  \int_{\xi_{1}+\xi_{2}=\xi}
\chi_{_{|\xi_{1}| \sim N_{1}}} \chi_{_{|\xi_{2}| \sim N_{2}}}
\chi_{_{| \cos{(\xi_{1},\xi_{2})}| \leq \theta}}
\chi_{_{[-\epsilon,\epsilon]}} (\tau - |\xi_{1}|^{2} -
|\xi_{2}|^{2}) \widehat{f}(\xi_{1}) \widehat{g}(\xi_{2}) d \xi_{1}
\end{array}
\label{Eqn:Def2Bil}
\end{equation}
then

\begin{equation}
\begin{array}{ll}
\overline{\lim} \frac{1}{2 \epsilon}  \|
B_{\epsilon}\|_{L_{\tau}^{2} L_{\xi}^{2}} & \lesssim \left(
\frac{N_{1}}{N_{2}} \right)^{\frac{1}{2}} \|  f \|_{L^{2}} \| g
\|_{L^{2}}
\end{array}
\label{Eqn:Bil}
\end{equation}
and

\begin{equation}
\begin{array}{ll}
\overline{\lim} \frac{1}{2 \epsilon}  \| B_{\epsilon, \theta}
\|_{L_{\tau}^{2} L_{\xi}^{2}} & \lesssim \theta^{\frac{1}{2}} \|  f
\|_{L^{2}} \| g \|_{L^{2}}.
\end{array}
\label{Eqn:AngleBil}
\end{equation}
The same conclusions hold if $-|\xi_{1}|^{2} -|\xi_{2}|^{2}$ is
substituted for $ |\xi_{1}|^{2} + |\xi_{2}|^{2}$ in
(\ref{Eqn:Def1Bil}) and (\ref{Eqn:Def2Bil}).
\end{prop}
Proposition \ref{prop:PropAlmCons}  shows that we can estimate the
variation of $\widetilde{E}(u(t))$ on an interval $J$  provided that
we can control the $X^{1,\frac{1}{2}+}$ norm of $Iu$. The next
proposition shows that such a control is possible as long as the
size of $J$ is, roughly speaking, bounded by one

\begin{prop}[Modified Local Well-Posedness in $X^{s,b}$ spaces \cite{cksttI2}]
There exists $1  < \epsilon \lesssim 1$ such that if  $\sup_{t \in
J} E(Iu(t)) \lesssim 1 $ and $|J| \leq \epsilon$ then
\begin{equation}
\begin{array}{ll}
\| \eta(t-a) I u \|_{X^{1,\frac{1}{2}+}} & \lesssim 1
\end{array}
\end{equation}
with $\eta$ bump function adapted to $[-\epsilon, \epsilon]$.
\label{Prop:ModLocWp}
\end{prop}
Finally, by choosing the optimal parameter $\theta= \frac{1}{N}$,
they estimated the variation of the almost conservation law
$\tilde{E}$ on an interval $J$ of size one

\begin{equation}
\begin{array}{ll}
| \sup_{t \in J} E(Iu(t)) - \widetilde{E}(u(a)) | & \lesssim
\frac{1}{N^{2-}}
\end{array}
\end{equation}
The variation is slower than that of the modified energy. Indeed
this $O \left( \frac{1}{N^{2-}} \right)$ increase is smaller than
the $O \left( \frac{1}{N^{\frac{3}{2}-}} \right)$ increase for the
variation of the modified energy \cite{cksttI1}.

\section{Summary of \cite{cgt}}
\label{sec:Morawetz}

In this section we recall two results from \cite{cgt} that we use in
the proof of Theorem \ref{Thm:Gwp13}. The first one shows that if we
the $L_{t}^{4} L_{x}^{4}$ norm of a solution to
(\ref{Eqn:Schroddefcubic} ) is small then we control several norms.
This result will be extensively used in establishing the
almost conservation law: see Proposition \ref{prop:VarNewAlmStr}.

\begin{prop}[Modified Local Well-Posedness \cite{cgt} ]
Let $u$ be a solution to ( \ref{Eqn:Schroddefcubic} ). Assume that
$(q,r)$ is admissible, i.e $(q,r) \in (2,\infty] \times [2,\infty) $
and $\frac{1}{q}+ \frac{1}{r}=\frac{1}{2}$. Assume also that

\begin{equation}
\begin{array}{ll}
\sup_{t \in J} E \left( I u (t) \right) & \leq 2
\end{array}
\label{Eqn:InducNrj}
\end{equation}
Then there exists $0 < \epsilon \ll 1 $ such that if

\begin{equation}
\begin{array}{ll}
\| I u \|_{L_{t}^{4}(J) L_{x}^{4}} & \leq \epsilon
\end{array}
\label{Eqn:Lt4Lx4Condition}
\end{equation}
then

\begin{equation}
\begin{array}{ll}
Z(J,u) & \lesssim 1
\end{array}
\label{Eqn:LocalWellPosedness}
\end{equation}
with

\begin{equation}
\begin{array}{ll}
Z(J,u) & := \sup_{(q,r) \, admissible} \| \langle D \rangle I u \|_{L_{t}^{q}(J)
L_{x}^{r}}
\end{array}
\end{equation}
\label{prop:LocalWellPosedness}
\end{prop}
The next result is a long-time estimate

\begin{prop}[Morawetz-type estimates \cite{cgt}, p9] Let $J$ be an interval and let $(J_{k})$ be a
partition of $J$. Let $u$ be the solution to
(\ref{Eqn:Schroddefcubic}) with data $u(0)=u_{0} \in H^{s}
(\mathbb{R}^{2})$. Then

\begin{equation}
\begin{array}{ll}
\| I u \|^{4}_{L_{t}^{4}(J) L_{x}^{4}} & \lesssim |J|^{\frac{1}{3}}
 \left( \sup_{t \in J} \| Iu(t) \|_{\dot{H}^{1}} \| I u \|^{3}_{L^{2}} + \| u_{0} \|^{4}_{L^{2}} +
\sum_{k} \frac{Z^{6}(J_{k},u)}{N^{1-}} \right)
\end{array}
\label{Eqn:Mor1}
\end{equation}
This inequality results from the two following estimates

\begin{equation}
\begin{array}{ll}
\| I u \|^{4}_{L_{t}^{4}(J) L_{x}^{4}} & \lesssim |J|^{\frac{1}{3}}
 \left( \sup_{t \in J} \| Iu(t) \|_{\dot{H}^{1}} \| I u \|^{3}_{L^{2}} + \| u_{0} \|^{4}_{L^{2}} +
Error(u,J) \right)
\end{array}
\label{Eqn:Mor2}
\end{equation}
and

\begin{equation}
\begin{array}{ll}
Error(u,J) & \lesssim \sum_{k} \frac{Z^{6}(J_{k},u)}{N^{1-}}
\end{array}
\label{Eqn:ErrMor}
\end{equation}
\label{prop:ImpMorEst}
\end{prop}

\section{Proof of global well-posedness in $H^{s}(\mathbb{R}^{2})$, $1 > s > \frac{1}{3} $}
\label{sec:PfGw13}

In this section we prove the global existence of
(\ref{Eqn:Schroddefcubic}) in  $H^{s}(\mathbb{R}^{2}) \times
H^{s-1}(\mathbb{R}^{2})$, $1 > s
> \frac{1}{3} $. Our proof relies on an intermediate result that we
prove in the next sections. More precisely we shall show the
following

\begin{prop}[Almost Conservation Law in $L_{t}^{q} L_{x}^{r}$ spaces]
Let $u$ be a solution of (\ref{Eqn:Schroddefcubic}). Assume that
(\ref{Eqn:InducNrj}) and (\ref{Eqn:Lt4Lx4Condition}) hold. Then

\begin{equation}
\begin{array}{ll}
\left| \sup_{t \in J} \widetilde{E} (u(t)) - \widetilde{E}(u(a))
\right| & \lesssim
\begin{array}{l}
\frac{1}{N^{2-}} +  \frac{ \theta^{\frac{1}{2}}}{N^{\frac{3}{2}-}} +
\frac{1}{\theta N^{3-}}
\end{array}
\end{array}
\label{Eqn:EstVarNewAlm}
\end{equation}
\label{prop:VarNewAlmStr}
\end{prop}
For the remainder of the section we show that Proposition
\ref{prop:VarNewAlmStr} implies Theorem \ref{Thm:Gwp13}.

Let $T>0$, $N=N(T) \gg 1 $ be a parameter to be chosen. There are
three steps to prove Theorem \ref{Thm:Gwp13}.

\begin{enumerate}

\item \textbf{Scaling}. We recall (see \cite{cksttI1})that there
exists $ C_{0}:=C_{0} \left( \| u_{0} \|_{H^{s}(\mathbb{R}^{2})}
\right)$ such that if $\lambda$ satisfies

\begin{equation}
\begin{array}{ll}
\lambda & =  C_{0} N^{\frac{1-s}{s}}
\end{array}
\label{Eqn:UpperBdLambda}
\end{equation}
then

\begin{equation}
\begin{array}{ll}
E \left( I u^{\lambda}(0) \right)  & \leq  \frac{1}{2}
\end{array}
\label{Eqn:InitTruncNrjEst}
\end{equation}
with $u^{\lambda}$ defined in (\ref{Eqn:Defulambda}).

\item \textbf{Bootstrap}. Let $F_{T}$ denote the following set

\begin{equation}
F_{T}= \left\{ T^{'} \in [0, \,T]:
\begin{array}{l}
\sup_{t \in [0, \, \lambda^{2}T^{'}]} E \left( I u^{\lambda}(t) \right) \leq 1 \,  \mathrm{and} \\
\| I u^{\lambda} \|^{4}_{L_{t}^{4}(J) L_{x}^{4}} \leq  2C_{mor}
|J|^{\frac{1}{3}} \max{ \left( \| u_{0} \|^{3}_{L^{2}}, \, \| u_{0}
\|^{4}_{L^{2}},1 \right) }, J \in \mathcal{P}_{N^{3--}}([0,
\lambda^{2} T^{'}])
\end{array}
\right\}
\end{equation}
with $\lambda$ defined in (\ref{Eqn:UpperBdLambda}) and $C_{mor}$
the constant determined by $\lesssim$ in (\ref{Eqn:Mor1}). We claim
that $F_{T}$ is the whole set $[0, \, T]$ for $N=N(T, \| u_{0}
\|_{H^{s}(\mathbb{R}^{2})}) \gg 1$ to be chosen later. Indeed

\begin{itemize}

\item $F_{T} \neq \emptyset$ since $0 \in F_{T}$ by
(\ref{Eqn:InitTruncNrjEst})

\item $F_{T}$ is closed by the dominated convergence theorem.

\item $F_{T}$ is open. Let $\widetilde{T^{'}} \in F_{T}$. Then by
continuity there exists $\delta>0$ such that for every $T^{'} \in
(\widetilde {T^{'}} - \delta, \widetilde{T^{'} } + \delta ) \cap [0,
T]$

\begin{equation}
\begin{array}{ll}
\sup_{t \in [0, \lambda^{2}T^{'}]} E(Iu^{\lambda}(t)) \leq 2
\end{array}
\label{Eqn:InductionNrj}
\end{equation}
and

\begin{equation}
\begin{array}{ll}
\| I u^{\lambda} \|^{4}_{L_{t}^{4}(J) L_{x}^{4}} & \leq 4C_{mor}
|J|^{\frac{1}{3}} \max{(\| u_{0} \|^{3}_{L^{2}}, \| u_{0}
\|^{4}_{L^{2}},1)}
\end{array}
\label{Eqn:AprioriLt4Lx4}
\end{equation}
for $J \in \mathcal{P}_{N^{3--}}([0, \lambda^{2} T^{'}])$.

Let $C_{fix}$ be the constant determined by $\lesssim$ in
(\ref{Eqn:FixedEst}).

Let $\theta=\frac{1}{N}$. Then by Proposition \ref{prop:FixedEst},
(\ref{Eqn:InitTruncNrjEst}) and the triangle inequality we have

\begin{equation}
\begin{array}{ll}
| \widetilde{E}(u^{\lambda}(0)) | & \leq \frac{1}{2} + \frac{C_{fix}}{N^{1-}} \\
& \leq \frac{5}{8}
\end{array}
\label{Eqn:EstTildeEinit}
\end{equation}
Then we divide each $J=[a_{j},b_{j}]$ into subintervals $J_{k}$, $k
\in \{ 1, \, ... l \} $ such that $ \| I u^{\lambda}
\|_{L_{t}^{4}(J_{k}) L_{x}^{4}}= \epsilon $, $k \in \{1, \, ... l-1
\}$ and $\| I u^{\lambda} \|_{L_{t}^{4}(J_{l}) L_{x}^{4}} \leq
\epsilon$ with $\epsilon$ defined in Proposition
\ref{prop:LocalWellPosedness}. By (\ref{Eqn:AprioriLt4Lx4}) we have

\begin{equation}
\begin{array}{ll}
l & \lesssim  N^{1--}
\end{array}
\label{Eqn:EstCardJk}
\end{equation}
By Proposition \ref{prop:LocalWellPosedness}, Proposition
\ref{prop:VarNewAlmStr}, (\ref{Eqn:EstTildeEinit}),
(\ref{Eqn:EstCardJk}) and by iteration we have

\begin{equation}
\begin{array}{ll}
\left| \sup_{t \in J } \widetilde{E}(u^{\lambda}(t))
- \widetilde{E}(u^{\lambda}(a_{j})) \right|   & \lesssim \frac{N}{N^{2-}} \\
& \lesssim \frac{1}{N^{1-}}
\end{array}
\label{Eqn:EstVarNewAlmJ}
\end{equation}
Now we iterate again to cover $[0, \lambda^{2} T^{'}]$. The number
of intervals $J$ is bounded by $\frac{\lambda^{2} T}{N^{3-}}$.
Therefore by this observation, (\ref{Eqn:EstTildeEinit}) and
(\ref{Eqn:EstVarNewAlmJ}) we have

\begin{equation}
\begin{array}{ll}
\left| \sup_{t \in [0,\lambda^{2} T^{'}] }
\widetilde{E}(u^{\lambda}(t)) \right| - \frac{5}{8} \lesssim
\frac{\lambda^{2}T}{N^{4-}}
\end{array}
\label{Eqn:TotVarNewAlm}
\end{equation}
By (\ref{Eqn:TotVarNewAlm}) and Proposition \ref{prop:FixedEst} we
have

\begin{equation}
\begin{array}{ll}
\left| \sup_{t \in [0, \lambda^{2} T] } E(Iu^{\lambda}(t)) \right| -
\frac{5}{8} \lesssim  \frac{\lambda^{2}T}{N^{4-}} + \frac{1}{N^{1-}}
\end{array}
\label{Eqn:EstTotNrj}
\end{equation}
Let $C_{tot}$ be the constant determined by $\lesssim$ in
(\ref{Eqn:EstTotNrj}). Since $s>\frac{1}{3}$ then for every $T>0$ we
can always choose $N=N(T) \gg 1$ such that $ C_{tot} \left(
\frac{\lambda^{2} T}{N^{4-}} + \frac{1}{N^{-}} \right) \leq
\frac{1}{8}$. Consequently $\sup_{t \in [0, \lambda^{2} T^{'} ]}
\widetilde{E}(u(t)) \lesssim 1$.

It remains to prove  $ \| I u^{\lambda} \|^{4}_{L_{t}^{4}(J)
L_{x}^{4}} \leq  4C_{mor} |J|^{\frac{1}{3}} \max{ \left( \| u_{0}
\|^{3}_{L^{2}}, \| u_{0} \|^{4}_{L^{2}} \right) }$, $J \in
\mathcal{P}_{N^{3--}}([0, \lambda^{2} T^{'}])$. We get from
Proposition \ref{prop:ImpMorEst}, (\ref{Eqn:EstCardJk}) and the
elementary inequality $\| I u(t) \|_{L^{2}} \leq \| u_{0} \|_{L^{2}}
$

\begin{equation}
\begin{array}{ll}
\| I u^{\lambda} \|^{4}_{L_{t}^{4}(J) L_{x}^{4}} - C_{mor}
|J|^{\frac{1}{3}} \left( \sqrt{2} \| u_{0} \|^{3}_{L^{2}} + \| u_{0}
\|^{4}_{L^{2}} \right) & \leq C_{mor}
\frac{|J|^{\frac{1}{3}}}{N^{+}}
\end{array}
\end{equation}
Hence

\begin{equation}
\begin{array}{ll}
\| I u^{\lambda} \|^{4}_{L_{t}^{4}(J) L_{x}^{4}} & \leq 2 C_{mor}
|J|^{\frac{1}{3}}  \max{ \left( \| u_{0} \|^{3}_{L^{2}}, \| u_{0}
\|^{4}_{L^{2}},1 \right) }
\end{array}
\end{equation}

\end{itemize}

\item \textbf{Accounting}. Following the $I$- method described
in \cite{cksttI1}

\begin{equation}
\begin{array}{ll}
\sup_{t \in [0, \, T]} E \left( I u(t) \right) & \lesssim
\lambda^{2}
\sup_{t \in [0, \, \lambda^{2} T]} E (I u^{\lambda}(t)) \\
& \lesssim \lambda^{2}
\end{array}
\label{Eqn:ComparNrjScale}
\end{equation}
By Proposition \ref{prop:HsnormNrj} we have global well-posedness of
the defocusing cubic Schr\"odinger equation in
$H^{s}(\mathbb{R}^{2})$, $1> s
> \frac{1}{3} $. Let $T \gg 1$. Then choosing $N=N(T) \gg 1$ such that

\begin{equation}
\begin{array}{ll}
\frac{0.9}{8} \leq C_{tot} \left(  \frac{\lambda^{2} T}{N^{4-}} +
\frac{1}{N^{1-}} \right) \leq \frac{1}{8}
\end{array}
\end{equation}
we have $N \sim T^{\frac{s}{6s-2}+} $. Plugging this value of $N$
into (\ref{Eqn:ComparNrjScale}) and using (\ref{Eqn:HsEIu}) we
obtain (\ref{Eqn:BoundHsuT}).

\end{enumerate}

\section{Proof of  Almost conservation law in $L_{t}^{q} L_{x}^{r}$ spaces}
\label{sec:AlmCon}

We modify an argument used in \cite{cksttI2}. Recall that the
derivative of $\tilde{E}$ is given by the following formula

\begin{equation}
\begin{array}{ll}
\partial_{t} \tilde{E} (u(t)) & = \Lambda_{4} \left( \mu + \sigma_{4}\alpha_{4},u(t)
\right) + \Lambda_{6} \left( \nu_{6}, u(t) \right).
\end{array}
\label{Eqn:DerivNewAlm}
\end{equation}
Let $J=[a, \, b]$ be an interval included in $[0,\infty)$ and let
$u$ be such that (\ref{Eqn:Schroddefcubic}), (\ref{Eqn:InducNrj}) and
(\ref{Eqn:Lt4Lx4Condition}) hold. The proof of the almost
conservation follows from (\ref{Eqn:DerivNewAlm}), the quadrilinear
estimate

\begin{equation}
\begin{array}{ll}
\left| \Lambda_{4,J} ( \mu + \sigma_{4} \alpha_{4} ,u ) \right| &
\lesssim \frac{\theta^{\frac{1}{2}}}{N^{\frac{3}{2}-}} +
\frac{1}{N^{2-}}
\end{array}
\label{Eqn:QuadEst}
\end{equation}
and the sextilinear estimate

\begin{equation}
\begin{array}{ll}
\left| \Lambda_{6,J} ( \sigma_{6},u ) \right| & \lesssim
\frac{1}{\theta N^{3-}}.
\end{array}
\label{Eqn:SextEst}
\end{equation}
Recall that $\alpha_{4}$, $\mu$ and $\sigma_{4}$ are defined in
(\ref{Eqn:Defalpha4}), (\ref{Eqn:Defmu}) and (\ref{Eqn:Defsigma4})
respectively.

\subsection{Proof of the quadrilinear estimate}

For convenience let $\nu_{4}$ denote the following multiplier

\begin{equation}
\begin{array}{ll}
\nu_{4} : =  \mu + \sigma_{4} \alpha_{4}
\end{array}
\end{equation}
so that by (\ref{Eqn:Defmu}) and (\ref{Eqn:Defalpha4}) we have

\begin{equation}
\begin{array}{ll}
\nu_{4}= \frac{i}{4} \left( |\xi_{1}|^{2} m^{2}(\xi_{1}) -
|\xi_{2}|^{2} m^{2}(\xi_{2}) + |\xi_{3}|^{2} m^{2}(\xi_{3})
-|\xi_{4}|^{2} m^{2}(\xi_{4}) \right) \chi_{L^{c} \cap \Gamma^{c}}.
\end{array}
\end{equation}

Notice that $\left| \Lambda_{4,J}  \right|$ is symmetric under
swapping $\xi_{1}$, $\xi_{2}$ with $\xi_{3}$, $\xi_{2}$
respectively. Therefore we may assume $|\xi_{1}| \geq  |\xi_{3}|$
and $|\xi_{2}| \geq |\xi_{4}|$. Notice that if we swap $\xi_{1}$,
$\xi_{3}$ with $\xi_{2}$, $\xi_{4}$ then $\left| \Lambda_{4,J}
\right|$ restricted to the set $\{(\xi_{1},...,\xi_{4}), \,
|\xi_{1}| \geq  |\xi_{3}|, \, |\xi_{2}| \geq |\xi_{4}| \}$ remains
invariant. Therefore we may also assume $|\xi_{3}| \geq |\xi_{4}|$.
Now we can restrict to $|\xi_{1}| \sim |\xi_{2}|$ since if not we
cannot have $| \cos{(\xi_{12},\xi_{14})} | \leq \theta$. Eventually
it suffices to prove

\begin{equation}
\begin{array}{ll}
\left| \Lambda_{4,J} ( \nu_{4} \chi_{_{\Sigma_{4}}}; u ) \right| &
\lesssim \frac{\theta^{\frac{1}{2}}}{N^{\frac{3}{2}-}} +
\frac{1}{N^{2-}}
\end{array}
\end{equation}
with

\begin{equation}
\begin{array}{ll}
\Sigma_{4} & = \left\{ (\xi_{1},...,\xi_{4}), \, |\xi_{1}| \gtrsim
N, \, |\xi_{1}| \sim |\xi_{2}|, \, |\xi_{1}| \geq |\xi_{3}| \geq
|\xi_{4}|, \, | \cos{(\xi_{12}, \xi_{14})} | \leq \theta \right\}.
\end{array}
\end{equation}
Then we need the following lemma

\begin{lem}

Let  $A \in \Omega_{k}^{1 \rightarrow 4}$, $k \in \{ 0,...,4 \}$.
Then

\begin{equation}
\begin{array}{ll}
\left| \Lambda_{4,J,A} ( \nu_{4}  \chi_{_{\Sigma_{4}}}; u) \right| &
\lesssim \left( \frac{\theta^{\frac{1}{2}}}{N^{\frac{3}{2}-}} +
\frac{1}{N^{2-}} \right) \| \langle D \rangle I u \|^{k}_{ L_{t}^{1}(d \mu^{l},J)
L_{x}^{2}} \| \langle D \rangle I (|u|^{2} u) \|^{4-k}_{L_{t}^{1}(J) L_{x}^{2}}
\end{array}
\label{Eqn:QuadIneq}
\end{equation}
\label{Lem:Quad}
\end{lem}
Let us postpone the proof of this lemma to later and let us assume
that it is true for the moment. Then we have

\begin{equation}
\begin{array}{ll}
\|\langle D \rangle I (|u|^{2} u) \|_{L_{t}^{1}(J) L_{x}^{2}} & \lesssim \| \langle D \rangle I
u \|_{L_{t}^{3}(J) L_{x}^{6}} \left(
\| P_{\ll N} u \|_{L_{t}^{3}(J) L_{x}^{6}} + \| P_{\gtrsim N} u \|_{L_{t}^{3}(J) L_{x}^{6}} \right)^{2} \\
& \lesssim \| \langle D \rangle I u \|_{L_{t}^{3}(J) L_{x}^{6}}  \left( \| \langle D \rangle I u
\|_{L_{t}^{3}(J) L_{x}^{6}} + \frac{\|\langle D \rangle I u  \|_{L_{t}^{3}(J)
L_{x}^{6}}}{N}
\right)^{2}  \\
& \lesssim Z^{3}(J,u)
\end{array}
\label{Eqn:Nlin}
\end{equation}
by the fractional Leibnitz rule and by H\"older inequality. Moreover

\begin{equation}
\begin{array}{ll}
\| \langle D \rangle I u \|_{L_{t}^{1} (J,d \mu^{l}) L_{x}^{2}} & \lesssim \| \langle D \rangle I u \|_{ L_{t}^{\infty}(J) L_{x}^{2}} \\
& \lesssim Z(J,u).
\end{array}
\label{Eqn:Lin}
\end{equation}
Therefore by Proposition \ref{prop:LocalWellPosedness}, Lemma
\ref{Lem:Quad}, (\ref{Eqn:Nlin}) and (\ref{Eqn:Lin}) we have

\begin{equation}
\begin{array}{ll}
\left| \Lambda_{4,J} ( \nu_{4} \chi_{_{\Sigma_{4}} }; u ) \right| &
\leq \sum_{k=0}^{4} \sum_{A \in \Omega_{k}^{1 \rightarrow 4}}
\left| \Lambda_{4,J,A} ( \nu_{4} \chi_{_{\Sigma_{4}}}; u) \right| \\
& \lesssim \left( \frac{\theta^{\frac{1}{2}}}{N^{\frac{3}{2}-}} +
\frac{1}{N^{2-}} \right) \sum_{k=0}^{4} \| \langle D \rangle I u \|^{k}_{
L_{t}^{1}(J,d \mu^{l}) L_{x}^{2}}
\| \langle D \rangle I (|u|^{2} u) \|^{4-k}_{L_{t}^{1}(J) L_{x}^{2}} \\
& \lesssim  \left( \frac{\theta^{\frac{1}{2}}}{N^{\frac{3}{2}-}} + \frac{1}{N^{2-}} \right)  \sum_{k=0}^{4} Z^{12-2k}(J,u) \\
& \lesssim \frac{\theta^{\frac{1}{2}}}{N^{\frac{3}{2}-}} +
\frac{1}{N^{2-}}.
\end{array}
\end{equation}
This proves the quadrilinear estimate (\ref{Eqn:QuadEst}).

\subsubsection{Proof of Lemma \ref{Lem:Quad}}

Given $k \in \{0,...,4 \}$ and $A \in \Omega^{1 \rightarrow 4}_{k}$
let $w_{j}$, $j \in \{1,...,4 \}$ denote the following functions

\begin{equation}
\begin{array}{rl}
w_{j}(t_{j}):= & \left\{
\begin{array}{l}
u(t_{j}), \, \mathrm{j \, odd}, \, \mathrm{j \in A} \\
|u|^{2} u(t_{j}), \, \mathrm{j \, odd}, \, \mathrm{j \notin A } \\
\overline{u}(t_{j}), \, \mathrm{j \, even}, \, \mathrm{j \in A }  \\
|u|^{2} \overline{u}(t_{j}), \, \mathrm{j \, even}, \, \mathrm{j
\notin A }
\end{array}
\right.
\end{array}
\end{equation}
and

\begin{equation}
\begin{array}{ll}
Q \left( t_{1},...,t_{4} \right) & := \int_{ \cap_{j=1}^{4}
[t_{j},b] } \int_{\xi_{1}+...+\xi_{4}=0} \nu_{4}
\chi_{_{\Sigma_{4}}} \prod_{1 \leq j \leq 4} e^{i
\epsilon(j)(t-t_{j})|\xi_{j}|^{2}} \widehat{w_{j}}(t_{j},\xi_{j}) dt
\end{array}
\end{equation}
with

\begin{equation}
\begin{array}{ll}
\epsilon(j) & = \left\{
\begin{array}{l}
1, \, j \, even \\
-1, \, j \, odd.
\end{array}
\right.
\end{array}
\label{Eqn:DefEpsilonj}
\end{equation}
We have

\begin{equation}
\begin{array}{ll}
\left| \Lambda_{4,J,4}(\nu_{4}\chi_{_{\Sigma_{4}}}  ; u)  \right| &
= \left| \int_{J} \int_{\xi_{1} +...+ \xi_{4}=0} \nu_{4}
\chi_{_{\Sigma_{4}}}
\begin{array}{l}
\left( \prod_{j \in A} \int_{a}^{t} e^{i
\epsilon(j)(t-t_{j})|\xi_{j}|^{2}} \widehat{w_{j}}(t_{j},\xi_{j})
d\mu_{l}(t_{j}) \right) \\
\left( \prod_{j \notin A} \int_{a}^{t} e^{i
\epsilon(j)(t-t_{j})|\xi_{j}|^{2}} \widehat{w_{j}}(t_{j},\xi_{j})
dt_{j} \right)
\end{array}
dt \right|
\end{array}
\end{equation}
and by Fubini

\begin{equation}
\begin{array}{ll}
\left| \Lambda_{4,J,A}(\nu_{4} \chi_{_{\Sigma_{4}}};u)  \right| & =
\left| \int_{J^{4}} Q(t_{1},...,t_{4},w_{1},...,w_{4}) \left(
\prod_{j \in A} d \mu_{l}(t_{j}) \right) \left( \prod_{j \notin A} d
t_{j} \right) \right|.
\end{array}
\label{Eqn:QuadExp}
\end{equation}
If we could prove

\begin{equation}
\begin{array}{ll}
| Q(t_{1},...,t_{4},w_{1},...,w_{4}) | & \lesssim \left(
\frac{\theta^{\frac{1}{2}}}{N^{\frac{3}{2}-}} + \frac{1}{N^{2-}}
\right) \prod_{j=1}^{4} \| \langle D \rangle I w_{j}(t_{j}) \|_{L_{x}^{2}}
\end{array}
\label{Eqn:BoundQ4}
\end{equation}
then (\ref{Eqn:QuadIneq}) would follow from (\ref{Eqn:QuadExp}) and
(\ref{Eqn:BoundQ4}).

We perform a Paley-Littlewood decomposition to prove
(\ref{Eqn:BoundQ4}). Let $X$ denote the left-hand side of
(\ref{Eqn:BoundQ4}) after decomposition. By Plancherel's theorem

\begin{equation}
\begin{array}{ll}
X  = & \left| \int_{_{ \tau_{0} + \tau^{'} + \tau^{''} =0}}
\widehat{\chi}_{_{\cap_{j=1}^{4} [t_{j},b]}}(\tau_{0})
\int_{_{\xi^{'} + \xi^{''}=0}}  \nu_{4} \chi_{_{\Sigma_{4}}}
\begin{array}{l}
\left( \int_{_{ \xi_{1} + \xi_{3} = \xi^{'}}}
\begin{array}{l}
\chi_{_{|\xi_{1}| \sim N_{1}}}  \chi_{_{|\xi_{3}| \sim N_{3}}}
\chi_{ _{|\cos{(\xi_{12},\xi_{14})}| \leq \theta}} \\ \\  \delta (
\tau^{'}  + |\xi_{1}|^{2} + |\xi_{3}|^{2} )
\widehat{w(t_{1})}(\xi_{1}) \widehat{w(t_{3})}(\xi_{3})
\end{array}
\right) \\
\left( \int_{_{\xi_{2} + \xi_{4} = \xi^{''}}}
\begin{array}{l}
\chi_{_{|\xi_{2}| \sim N_{2}}} \chi_{_{|\xi_{4}| \sim N_{4}}} \delta
\left( \tau^{''} - |\xi_{2}|^{2} - |\xi_{4}|^{2} \right) \\
\\ \widehat{w(t_{2})}(\xi_{2}) \widehat{w(t_{4})}(\xi_{4})
\end{array}
\right)
\end{array}
\right|
\end{array}
\end{equation}
and we want to prove

\begin{equation}
\begin{array}{ll}
X & \lesssim   N_{1}^{--} N_{4}^{+} \left(
\frac{\theta^{\frac{1}{2}}}{N^{\frac{3}{2}-}} + \frac{1}{N^{2-}}
\right) \prod_{j=1}^{4} \|\langle D \rangle I w_{j}(t_{j}) \|_{L_{x}^{2}}.
\end{array}
\end{equation}
Since the $L^{2}$- norm only depends on the magnitude of the Fourier
transform we may assume that $\widehat{w(t_{j})} \geq 0$, $j \in
\{1,...,4 \}$. There are two cases:

\begin{itemize}

\item  \textbf{Case} $\mathbf{1}$: $N_{3} \gg N_{4}$. Recall (see \cite{cksttI2}) that

\begin{equation}
\begin{array}{ll}
\left| \cos{(\xi_{1},\xi_{3})} \right| & \lesssim \theta_{0} + \frac{N_{4}}{N_{3}} \\
\end{array}
\label{Eqn:BoundCos}
\end{equation}
and

\begin{equation}
\begin{array}{ll}
\| \nu_{4} \chi_{_{\Sigma_{4}}} \|_{L^{\infty}} & \leq m^{2}(N_{1})
N_{1} N_{3} \theta + m^{2}(N_{3}) N_{3}^{2}
\end{array}
\end{equation}
There are two subcases

\begin{itemize}

\item \textbf{Case} $\mathbf{1.a}$: $\theta \gtrsim \frac{N_{4}}{N_{3}}$

We have

\begin{equation}
\begin{array}{ll}
\left| \widehat{\chi}_{_{\cap_{j=1}^{4} [t_{j},b]}}(\tau_{0})
\right| & \lesssim <\tau_{0}>
\end{array}
\label{Eqn:BdChi}
\end{equation}
We introduce the logarithmic weight

\begin{equation}
\begin{array}{ll}
q(\tau) & := 1+ \log^{2}(<\tau>).
\end{array}
\label{Eqn:DefWeight}
\end{equation}
Notice that $q(\tau^{'} + \tau^{''}) \lesssim q(\tau^{'}) +
q(\tau^{''})$. Let

\begin{equation}
\begin{array}{ll}
\widetilde{B_{_{\epsilon,1,3}}} (\tau^{'},\xi^{'}) & : =
\int_{\xi_{1} + \xi_{3}= \xi^{'}} \chi_{_{|\xi_{1}| \sim N_{1}} }
\chi_{ _{|\xi_{3}| \sim N_{3}} } \chi_{[-\epsilon, \epsilon]}
(\tau^{'} - |\xi_{1}|^{2} - |\xi_{3}|^{2}) \chi_{|
\cos{(\xi_{1},\xi_{3})|} \lesssim \theta_{0}}
\widehat{w(t_{1})}(\xi_{1}) \widehat{w(t_{3})}(\xi_{3})
\end{array}
\label{Eqn:DefBepsilon13Angle}
\end{equation}
and

\begin{equation}
\begin{array}{ll}
\widetilde{B_{_{\epsilon,2,4}}} (\tau^{''},\xi^{''}) & : =
\int_{\xi_{2} + \xi_{4}= \xi^{''}} \chi_{_{|\xi_{2}| \sim N_{2}}}
\chi_{_{|\xi_{2}| \sim N_{4}}} \chi_{[-\epsilon, \epsilon]}
(\tau^{''} - |\xi_{2}|^{2} - |\xi_{4}|^{2})
\widehat{w(t_{2})}(\xi_{2}) \widehat{w(t_{4})}(\xi_{4}).
\end{array}
\label{Eqn:DefBepsilon24}
\end{equation}
Then by Hausdorff-Young, (\ref{Eqn:Bil}), (\ref{Eqn:AngleBil}) and
(\ref{Eqn:BdChi})

\begin{equation}
\begin{array}{ll}
X & \lesssim q^{2}(N_{1}) \| \nu_{4} \chi_{_{\Sigma_{4}}}
\|_{L^{\infty}} \lim_{\epsilon \rightarrow 0} \frac{1}{(2
\epsilon)^{2}} \int_{\mathbb{R}} \frac{1}{<\tau_{0}> q(\tau_{0})}
\left ( \widetilde{B_{_{\epsilon,1,3}}}
* \widetilde{B_{_{\epsilon,2,4}}} (\tau_{0},0) \right) \, d \tau_{0} \\
& \lesssim q^{2}(N_{1}) \| \nu_{4} \chi_{_{\Sigma_{4}}}
\|_{L^{\infty}} \overline{\lim}  \frac{1}{(2 \epsilon)^{2}}
\| \widetilde{B_{_{\epsilon,1,3}}} * \widetilde{B_{_{\epsilon,2,4}}} \|_{L^{\infty}_{\tau_{0}} L^{\infty}_{\xi_{0}}} \\
& \lesssim q^{2}(N_{1}) \| \nu_{4} \chi_{_{\Sigma_{4}}}
\|_{L^{\infty}} \overline{\lim}  \frac{1}{(2 \epsilon)^{2}}
\| B_{_{\epsilon,1,3}} B_{_{\epsilon,2,4}} \|_{L^{1}_{t} L^{1}_{x}} \\
& \lesssim q^{2}(N_{1}) \| \nu_{4} \chi_{_{\Sigma_{4}}}
\|_{L^{\infty}} \overline{\lim}  \frac{1}{(2 \epsilon)^{2}}  \|
B_{_{\epsilon,1,3}} \|_{L_{t}^{2} L_{x}^{2}}
\|  B_{_{\epsilon,2,4}}  \|_{L_{t}^{2} L_{x}^{2}} \\
& \lesssim q^{2}(N_{1}) \frac{ m^{2}(N_{1}) N_{1} N_{3} \theta +
m^{2}(N_{3}) N_{3}^{2}} {<N_{1}> m(N_{1})...<N_{4}> m(N_{4})}
\theta^{\frac{1}{2}}
\left( \frac{N_{4}}{N_{2}} \right)^{\frac{1}{2}}  \prod_{j=1}^{4} \| \langle D \rangle I w_{j}(t_{j}) \|_{L_{x}^{2}} \\
& \lesssim N_{1}^{--} N_{4}^{+}
\frac{\theta^{\frac{1}{2}}}{N^{\frac{3}{2}-}} \prod_{j=1}^{4} \| \langle D \rangle
I w_{j}(t_{j}) \|_{L_{x}^{2}}
\end{array}
\label{Eqn:DvtQuad}
\end{equation}

\item \textbf{Case} $\mathbf{1.b}$: $\theta \ll \frac{N_{4}}{N_{3}}$

In this case $|\cos{(\xi_{1},\xi_{3})}| \lesssim
\frac{N_{4}}{N_{3}}$ and

\begin{equation}
\begin{array}{ll}
X & \lesssim q^{2}(N_{1}) \| \nu_{4} \chi_{_{\Sigma_{4}}}
\|_{L^{\infty}} \overline{\lim} \frac{1}{(2 \epsilon)^{2}}
\| B_{_{\epsilon,1,3}} B_{_{\epsilon,2,4}} \|_{L^{1}_{t} L^{1}_{x}} \\
& \lesssim q^{2}(N_{1}) \frac{ m^{2}(N_{1}) N_{1} N_{3} \theta +
m^{2}(N_{3}) N_{3}^{2}} {<N_{1}> m(N_{1})...<N_{4}> m(N_{4})} \left(
\frac{N_{4}}{N_{3}} \right)^{\frac{1}{2}}
\left( \frac{N_{4}}{N_{2}} \right)^{\frac{1}{2}}  \prod_{j=1}^{4} \| \langle D \rangle I w_{j}(t_{j}) \|_{L_{x}^{2}} \\
& \lesssim N_{1}^{--} N_{4}^{+} \left(
\frac{\theta^{\frac{1}{2}}}{N^{\frac{3}{2}-}} + \frac{1}{N^{2-}}
\right) \prod_{j=1}^{6} \| \langle D \rangle I w_{j}(t_{j}) \|_{L_{x}^{2}}
\end{array}
\end{equation}
with

\begin{equation}
\begin{array}{ll}
\widetilde{B_{_{\epsilon,1,3}}} (\tau^{'},\xi^{'}) & : = q(\tau^{'})
\int_{\xi_{1} + \xi_{3}= \xi^{'}} \chi_{_{|\xi_{1}| \sim N_{1}} }
\chi_{ _{|\xi_{3}| \sim N_{3}} } \chi_{[-\epsilon, \epsilon]}
(\tau^{'} - |\xi_{1}|^{2} - |\xi_{3}|^{2}) \chi_{|
\cos{(\xi_{1},\xi_{3})|} \lesssim \frac{N_{4}}{N_{3}}}
\widehat{w(t_{1})}(\xi_{1}) \widehat{w(t_{3})}(\xi_{3})
\end{array}
\label{Eqn:DefBepsilon13AngleOther}
\end{equation}
and $\widetilde{B_{_{\epsilon,2,4}}}$ defined in
(\ref{Eqn:DefBepsilon24}).

\end{itemize}

\item\textbf{Case} $\mathbf{2}$: $N_{3} \sim N_{4}$. Recall (see \cite{cksttI2}) that

\begin{equation}
\begin{array}{ll}
\| \nu_{4} \chi_{_{\Sigma_{4}}} \|_{L^{\infty}} & \lesssim
m^{2}(N_{1}) N_{1} N_{3}.
\end{array}
\end{equation}
We have

\begin{equation}
\begin{array}{ll}
X & \lesssim q^{2}(N_{1}) \| \nu_{4} \chi_{_{\Sigma_{4}}}
\|_{L^{\infty}} \overline{\lim} \frac{1}{(2 \epsilon)^{2}}
\| B_{_{\epsilon,1,3}} B_{_{\epsilon,2,4}} \|_{L^{1}_{t} L^{1}_{x}} \\
& \lesssim q^{2}(N_{1}) \frac{ m^{2}(N_{1}) N_{1} N_{3} } {<N_{1}>
m(N_{1})...<N_{4}> m(N_{4})} \left( \frac{N_{3}}{N_{1}}
\right)^{\frac{1}{2}}
\left( \frac{N_{4}}{N_{2}} \right)^{\frac{1}{2}}  \prod_{j=1}^{4} \| \langle D \rangle I w_{j}(t_{j}) \|_{L_{x}^{2}} \\
& \lesssim N_{1}^{--} N_{4}^{+} \frac{1}{N^{2-}} \prod_{j=1}^{4} \|
\langle D \rangle I w_{j}(t_{j}) \|_{L_{x}^{2}}
\end{array}
\end{equation}
with

\begin{equation}
\begin{array}{ll}
\widetilde{B_{_{\epsilon,1,3}}} (\tau^{'},\xi^{'}) & : = q(\tau^{'})
\int_{\xi_{1} + \xi_{3}= \xi^{'}} \chi_{_{|\xi_{1}| \sim N_{1}} }
\chi_{ _{|\xi_{3}| \sim N_{3}} } \chi_{[-\epsilon, \epsilon]}
(\tau^{'} + |\xi_{1}|^{2} + |\xi_{3}|^{2})
\widehat{w_{1}(t_{1})}(\xi_{1}) \widehat{w_{3}(t_{3})}(\xi_{3})
\end{array}
\label{Eqn:DefotherBepsilon13}
\end{equation}
and $\widetilde{B_{_{\epsilon,2,4}}}$ defined in
(\ref{Eqn:DefBepsilon24}).

\end{itemize}

\subsection{Proof of the sextilinear estimate}

Notice that $\nu_{6}=0$ if $\max{(|\xi_{1}|,...,|\xi_{6}|)} \ll N$.
Let $|\xi_{1^{*}}| \geq ... \geq |\xi_{6^{*}}|$ be the six
amplitudes in order. The convolution constraint
$\xi_{1}+...+\xi_{6}=0$ imposes $|\xi_{1^{*}}| \sim |\xi_{2^{*}}|$.
It suffices to prove

\begin{equation}
\begin{array}{ll}
\left| \Lambda_{6,J}(\nu_{6} \chi_{_{\Sigma_{6}}},u) \right|  &
\lesssim \frac{1}{\theta N^{3-}}
\end{array}
\end{equation}
with

\begin{equation}
\begin{array}{ll}
\Sigma_{6} & = \left\{ (\xi_{1},...,\xi_{6}), \, |\xi_{1^{*}}|
\gtrsim N, \, |\xi_{1^{*}}| \sim |\xi_{2^{*}}| \right\}.
\end{array}
\end{equation}
We will prove the following lemma

\begin{lem}
Let  $A \in \Omega_{k}^{1 \rightarrow 6}$, $k \in \{ 0,...,6 \}$.
Then

\begin{equation}
\begin{array}{ll}
\left| \Lambda_{6,J,A} (\nu_{6} \chi_{_{\Sigma_{6}}} ; u) \right| &
\lesssim \frac{1}{\theta N^{3-}} \| \langle D \rangle I u \|^{k}_{ L_{t}^{1}(d
\mu^{l},J) L_{x}^{2}} \| \langle D \rangle I (|u|^{2} u) \|^{6-k}_{L_{t}^{1}(J)
L_{x}^{2}}
\end{array}
\label{Eqn:SexIneq}
\end{equation}
\label{Lem:Sex}
\end{lem}
Assuming that it is true then by (\ref{Eqn:Nlin}), (\ref{Eqn:Lin})
and Proposition \ref{prop:LocalWellPosedness} we have

\begin{equation}
\begin{array}{ll}
\left| \Lambda_{6,J} ( \nu_{6} \chi_{\Sigma_{6}}; u ) \right| & \leq
\sum_{k=0}^{6} \sum_{A \in \Omega_{k}^{1 \rightarrow 6}}
\left| \Lambda_{6,J,A} (\nu_{6} \chi_{_{\Sigma_{6}}}; u) \right| \\
& \lesssim \frac{1}{\theta N^{3-}} \sum_{k=0}^{6} \| \langle D \rangle I u
\|^{k}_{ L_{t}^{1}(J,d \mu^{l}) L_{x}^{2}}
\| \langle D \rangle I (|u|^{2} u) \|^{6-k}_{L_{t}^{1}(J) L_{x}^{2}} \\
& \lesssim \frac{1}{\theta N^{3-}} \sum_{k=0}^{6} Z^{18-2k}(J,u) \\
& \lesssim \frac{1}{N^{3-}}
\end{array}
\end{equation}
which proves the sextilinear estimates.

\subsubsection{Proof of Lemma \ref{Lem:Sex}}

Given $k \in \{0,...,6 \}$ and $A \in \Omega^{1 \rightarrow 6}_{k}$
let $w_{j}$, $j \in \{1,...,6 \}$ denote the following functions

\begin{equation}
\begin{array}{rl}
w_{j}(t_{j}):= & \left\{
\begin{array}{l}
u(t_{j}), \, \mathrm{j \, odd}, \, \mathrm{j \in A} \\
|u|^{2} u(t_{j}), \, \mathrm{j \, even}, \, \mathrm{j \notin A } \\
\overline{u}(t_{j}), \, \mathrm{j \, even}, \, \mathrm{j \in A }  \\
|u|^{2} \overline{u}(t_{j}), \, \mathrm{j \, even}, \, \mathrm{j
\notin A }
\end{array}
\right.
\end{array}
\end{equation}
and

\begin{equation}
\begin{array}{ll}
Q \left( t_{1},...,t_{6};w_{1},...,w_{6} \right) & := \int_{
\cap_{j=1}^{6} [t_{j},b] } \int_{\xi_{1}+...+\xi_{6}=0} \nu_{6}
\chi_{_{\Sigma_{6}}} \prod_{1 \leq j \leq 6} e^{i
\epsilon(j)(t-t_{j})|\xi_{j}|^{2}} \widehat{w_{j}}(t_{j},\xi_{j}) dt
\end{array}
\end{equation}
with $\epsilon(j)$ defined in (\ref{Eqn:DefEpsilonj}). We have

\begin{equation}
\begin{array}{ll}
\left| \Lambda_{6,J,A}(\nu_{6}; u)  \right| & = \left| \int_{J}
\int_{\xi_{1} +...+ \xi_{6}=0} \nu_{6}  \left(
\begin{array}{l}
\left( \prod_{j \in A} \int_{a}^{t} e^{i
\epsilon(j)(t-t_{j})|\xi_{j}|^{2}}
\widehat{w_{j}}(t_{j},\xi_{j}) d \mu_{l}(t_{j}) \right) \\
\left( \prod_{j \notin A} \int_{a}^{t} e^{i
\epsilon(j)(t-t_{j})|\xi_{j}|^{2}} \widehat{w_{j}}(t_{j},\xi_{j})
dt_{j} \right)
\end{array}
\right) dt \right|
\end{array}
\end{equation}
and by Fubini

\begin{equation}
\begin{array}{ll}
\left| \Lambda_{6,J,A}(\nu_{6} \chi_{_{\Sigma_{6}}};u)  \right| & =
\left| \int_{J^{6}} Q(t_{1},...,t_{6},w_{1},...,w_{6}) \left(
\prod_{j \in A} d \mu_{l}(t_{j}) \right) \left( \prod_{j \notin A} d
t_{j} \right) \right|.
\end{array}
\label{Eqn:SextExp}
\end{equation}
If we could prove

\begin{equation}
\begin{array}{ll}
| Q(t_{1},...,t_{6},w_{1},...,w_{6}) | & \lesssim \frac{1}{\theta
N^{3-}} \prod_{j=1}^{6} \| \langle D \rangle I w_{j}(t_{j}) \|_{L_{x}^{2}}
\end{array}
\label{Eqn:BoundQ6}
\end{equation}
then (\ref{Eqn:SexIneq}) would follow from (\ref{Eqn:SextExp}) and
(\ref{Eqn:BoundQ6}). It remains to show (\ref{Eqn:BoundQ6}). By
decomposition we may assume $ \widehat{w_{j}(t_{j})} \geq 0$.

We perform a Paley-Littlewood decomposition to prove
(\ref{Eqn:BoundQ6}). Let $X$ be the left-hand side of
(\ref{Eqn:BoundQ6}). By Plancherel we have

\begin{equation}
\begin{array}{ll}
X = \left|  \int_{_ {*_{\tau}=0}} \widehat{\chi_{_{\cap_{j=1}^{6}
[t_{j},b]}}}(\tau_{0}) \int_{_{*_{\xi}=0}} \nu_{6} \left[
\begin{array}{l}
\left( \int_{_{ \xi_{1^{*}} + \xi_{3^{*}} = \xi^{'}}} \left[
\begin{array}{l}
\chi_{_{|\xi_{1^{*}}| \sim N_{1^{*}}}}  \chi_{_{|\xi_{3^{*}}| \sim N_{3^{*}}}}  \\
\delta ( \tau^{'} \pm |\xi_{1^{*}}|^{2} \pm |\xi_{3^{*}}|^{2} )
\widehat{w_{1}}(t_{1^{*}},\xi_{1^{*}})
\widehat{w_{3}}(t_{3^{*}},\xi_{3^{*}})
\end{array}
\right]
\right) \\
\left( \int_{_{\xi_{2^{*}} + \xi_{4^{*}} = \xi^{''}}} \left[
\begin{array}{l}
\chi_{_{|\xi_{2^{*}}| \sim N_{2^{*}}}} \chi_{_{|\xi_{4^{*}}| \sim N_{4^{*}}}}  \\
\delta \left( \tau^{''} \pm |\xi_{2^{*}}|^{2} \pm |\xi_{4^{*}}|^{2}
\right) \widehat{w_{2}}(t_{2^{*}},\xi_{2^{*}})
\widehat{w_{4}}(t_{4^{*}},\xi_{4^{*}})
\end{array}
\right]
\right) \\
\chi_{_{|\xi_{5^{*}}| \sim N_{5^{*}}}} \delta \left( \tau_{5^{*}} \pm |\xi_{5^{*}}|^{2} \right) \widehat{w_{5^{*}}}(t_{5^{*}},\xi_{5^{*}}) \\
\chi_{_{|\xi_{6^{*}}| \sim N_{6^{*}}}}  \delta \left( \tau_{6^{*}}
\pm |\xi_{6^{*}}|^{2} \right)
\widehat{w_{6^{*}}}(t_{6^{*}},\xi_{6^{*}})
\end{array}
\right] \right|
\\
\end{array}
\label{Eqn:SextPlanch}
\end{equation}
where $N_{1^{*}} \geq ... \geq N_{6^{*}}$ are the dyadic numbers in
order, $1^{*}$,..., $6^{*}$ are the corresponding subscripts,
$*_{\tau}:=\tau_{0}+ \tau^{'}+ \tau^{''}+\tau_{5^{*}}
+\tau_{6^{*}}$,  $*_{\xi}:=\xi^{'}+ \xi^{''}+ \xi_{5^{*}} +
\xi_{6^{*}}$, $\pm |\xi_{j}|^{2}$ denotes $+ |\xi_{j}|^{2}$ if $j$
is odd and $- |\xi_{j}|^{2}$ if $j$ is even. We would like to prove

\begin{equation}
\begin{array}{ll}
X & \lesssim \frac{N_{1^{*}}^{--} N_{6^{*}}^{+}}{\theta N^{3-}}
\prod_{j=1}^{6} \| \langle D \rangle I w_{j}(t_{j}) \|_{L_{x}^{2}}
\end{array}
\end{equation}
Again we can assume the $\widehat{w_{j}}(t_{j}) \geq 0$. Notice that

\begin{equation}
\begin{array}{ll}
\widehat{\chi_{_{\cap_{j=1}^{6} [t_{j},b]}}}(\tau_{0}) & \lesssim
\langle \tau_{0}\rangle^{-1}.
\end{array}
\label{Eqn:Charac6}
\end{equation}
Recall also (see \cite{cksttI2}) that $|
\nu_{4}(\xi_{1},...,\xi_{4}) | \lesssim \frac{\min{\left(
m(N_{1}),...,m(N_{4}) \right)^{2}}}{\theta}$. Therefore

\begin{equation}
\begin{array}{ll}
|\nu_{6}| & \lesssim
\sum_{k=1}^{4} \frac{ \left| \min{ ( m^{2}(\xi_{1}),...,m^{2}(\xi_{k}+..+\xi_{k+2})...,m^{2}(\xi_{6}) )} \right|}{\theta} \\
& \lesssim \frac{m^{2}(N_{4^{*}})}{\theta}.
\end{array}
\label{Eqn:EstMult6}
\end{equation}
Before continuing we define  $M_{\epsilon,j}$ and $P_{N_{j}}$  such
that

\begin{equation}
\begin{array}{ll}
\widetilde{M_{\epsilon,j}}(\tau_{j},\xi_{j}) :=
\chi_{[-\epsilon,\epsilon]} (\tau_{j} \pm |\xi_{j}|^{2})
\chi_{_{|\xi_{j}| \sim N_{j}}} \widehat{w_{j}}(t_{j},\xi_{j})
\end{array}
\end{equation}
and

\begin{equation}
\begin{array}{ll}
\widehat{P_{N_{j}}(f)}(\xi_{j}) := \chi_{|\xi_{j}| \sim N_{j}}
\hat{f}(\xi_{j})
\end{array}
\end{equation}
for $j \in [1,...,6 ]$. Also let $B_{_{\epsilon,k^{*},l^{*}}}$ be
such that

\begin{equation}
\begin{array}{l}
\widetilde{B_{_{\epsilon,k^{*},l^{*}}}}(\tau^{'},\epsilon^{'}) : =
\int_{\xi_{k^{*}} + \xi_{l^{*}}= \xi^{'}}
\begin{array}{l}
\chi_{_{|\xi_{k^{*}}| \sim N_{k^{*}}} } \chi_{ _{|\xi_{l^{*}}| \sim
N_{l^{*}}} } \chi_{[-\epsilon, \epsilon]} (\tau^{'} \pm |\xi_{k^{*}}|^{2} \pm |\xi_{l^{*}}|^{2}) \\
\widehat{w_{k^{*}}(t_{k^{*}})}(\xi_{k^{*}})
\widehat{w_{l^{*}}(t_{l^{*}})}(\xi_{l^{*}}).
\end{array}
\end{array}
\label{Eqn:DefotherBepsilonkstarlstar}
\end{equation}
Then we prove the following claim.

\underline{Claim}: If  $N_{k^{*}} \leq N_{l^{*}}$ then

\begin{equation}
\begin{array}{ll}
\overline{\lim} \frac{1}{2 \epsilon} \| B_{_{\epsilon,k^{*},l^{*}}}
\|_{L_{t}^{2} L_{x}^{2}} & \lesssim \left(
\frac{N_{k^{*}}}{N_{l^{*}}} \right)^{\frac{1}{2}} \|
w_{k^{*}}(t_{k^{*}}) \|_{L^{2}} \| w_{l^{*}}(t_{l^{*}}) \|_{L^{2}}.
\end{array}
\end{equation}

\begin{proof}
If $k^{*}$ and $l^{*}$ are of the same parity then then the claim
follows from Proposition \ref{Eqn:AngleBil}. It remains to study the
case where $k^{*}$ and $l^{*}$ are of different parity. Let
$B_{+,k,l}$, $B_{-,k,l}$ be such that

\begin{equation}
\begin{array}{l}
\widetilde{B_{_{+,\epsilon,k^{*},l^{*}}}}(\tau^{'},\epsilon^{'}) : =
\int_{\xi_{k^{*}} + \xi_{l^{*}}= \xi^{'}}
\begin{array}{l}
\chi_{_{|\xi_{k^{*}}| \sim N_{k^{*}}} } \chi_{ _{|\xi_{l^{*}}| \sim
N_{l^{*}}} } \chi_{[-\epsilon, \epsilon]} (\tau^{'} +  |\xi_{k^{*}}|^{2} + |\xi_{l^{*}}|^{2}) \\
\widehat{w_{k^{*}}(t_{k^{*}})}(\xi_{k^{*}})
\widehat{\overline{w_{l^{*}}}(t_{l^{*}})}(\xi_{l^{*}})
\end{array}
\end{array}
\end{equation}
and
\begin{equation}
\begin{array}{l}
\widetilde{B_{_{-,\epsilon,k^{*},l^{*}}}}(\tau^{'},\epsilon^{'}) : =
\int_{\xi_{k^{*}} + \xi_{l^{*}}= \xi^{'}}
\begin{array}{l}
\chi_{_{|\xi_{k^{*}}| \sim N_{k^{*}}} } \chi_{ _{|\xi_{l^{*}}| \sim
N_{l^{*}}} } \chi_{[-\epsilon, \epsilon]} (\tau^{'} -  |\xi_{k^{*}}|^{2} - |\xi_{l^{*}}|^{2}) \\
\widehat{w_{k^{*}}(t_{k^{*}})}(\xi_{k^{*}})
\widehat{\overline{w_{l^{*}}}(t_{l^{*}})}(\xi_{l^{*}}).
\end{array}
\end{array}
\end{equation}
Observe that

\begin{equation}
\begin{array}{ll}
\overline{\lim} \frac{1}{2 \epsilon} \| B_{_{\epsilon,k^{*},l^{*}}}
\|_{L_{t}^{2} L_{x}^{2}} & = \| P_{N_{k^{*}}} \left( e^{it
\epsilon(k^{*}) \triangle} w_{k^{*}}(t_{k^{*}})
\right) P_{N_{l^{*}}} \left(  e^{it \epsilon(l^{*})  \triangle } w_{l^{*}}(t_{l^{*}}) \right) \|_{L_{t}^{2} L_{x}^{2}} \\
& = \| P_{N_{k^{*}}} \left( e^{it \epsilon(k^{*}) \triangle}
w_{k^{*}}(t_{k^{*}}) \right) \overline{P_{N_{l^{*}}} \left( e^{ it
\epsilon(l^{*}) \triangle}
w_{l^{*}}(t_{l^{*}}) \right)} \|_{L_{t}^{2} L_{x}^{2}} \\
& \lesssim \overline{\lim} \frac{1}{2 \epsilon} \left( \|
\widetilde{B_{_{+,\epsilon,k^{*},l^{*}}}} \|_{L^{2}_{\tau} L^{2}_{\epsilon}} + \| \widetilde{B_{_{-,\epsilon,k^{*},l^{*}}}} \|_{L^{2}_{\tau} L^{2}_{\epsilon}} \right) \\
& \lesssim \left( \frac{N_{k^{*}}}{N_{l^{*}}} \right)^{\frac{1}{2}}
\| w_{k^{*}}(t_{k^{*}}) \|_{L^{2}} \|
\overline{w_{l^{*}}}(t_{l^{*}}) \|_{L^{2}} \\
& \lesssim \left( \frac{N_{k^{*}}}{N_{l^{*}}} \right)^{\frac{1}{2}}
\| w_{k^{*}}(t_{k^{*}}) \|_{L^{2}} \| w_{l^{*}}(t_{l^{*}}).
\|_{L^{2}}
\end{array}
\end{equation}
This ends the proof of the claim.
\end{proof}
Observe also that

\begin{equation}
\begin{array}{ll}
\overline{\lim} \frac{1}{2 \epsilon} \| M_{\epsilon,j}
\|_{L_{t}^{\infty} L_{x}^{\infty}} & \lesssim
\| e^{it \pm \triangle}  \left( P_{N_{j}} w_{j}(t_{j}) \right) \|_{L_{t}^{\infty} L_{x}^{\infty}} \\
& \lesssim \frac{N_{j}^{\frac{3}{2}}}  { m(N_{j}) <N_{j}> } \| \langle D \rangle I
w_{j}(t_{j})\|_{L_{x}^{2}}
\end{array}
\label{Eqn:MepsilonIneq}
\end{equation}
by Plancherel and Bernstein inequalities.

By (\ref{Eqn:Charac6}), (\ref{Eqn:EstMult6}),
(\ref{Eqn:MepsilonIneq}), the claim and Haussdorf-Young  we have

\begin{equation}
\begin{array}{ll}
X  & \lesssim \frac{ m^{2}(N_{4^{*}})}{\theta} q^{4}(N_{1^{*}})
\overline{\lim} \frac{1}{(2 \epsilon)^{4}} \|
\widetilde{B_{\epsilon,1^{*},3^{*}}} *
\widetilde{B_{\epsilon,2^{*},4^{*}}}
* \widetilde{M_{\epsilon,5^{*}}} * \widetilde{M_{\epsilon,6^{*}}} \|_
{L^{\infty}_{\tau} L_{\xi}^{\infty}} \\
& \lesssim  \frac{m^{2}(N_{4^{*}})}{\theta} q^{4}(N_{1^{*}})
\overline{\lim} \frac{1}{(2 \epsilon)^{4}}
\|B_{\epsilon,1^{*},3^{*}} B_{\epsilon,2^{*},4^{*}}
M_{\epsilon,5^{*}}
M_{\epsilon,6^{*}} \|_{L_{t}^{1} L_{x}^{1}} \\
& \lesssim q^{4}(N_{1^{*}}) \frac{m^{2}(N_{4^{*}})}{\theta}
\overline{\lim} \frac{1}{(2 \epsilon)^{4}} \|
B_{\epsilon,1^{*},3^{*}} \|_{L_{t}^{2} L_{x}^{2}} \|
B_{\epsilon,2^{*},4^{*}} \|_{L_{t}^{2} L_{x}^{2}}
\| M_{\epsilon,5^{*}} \|_{L_{t}^{\infty} L_{x}^{\infty}} \| M_{\epsilon,6^{*}} \|_{L_{t}^{\infty} L_{x}^{\infty}} \\
&  \lesssim q^{4}(N_{1^{*}}) \frac{m^{2}(N_{4^{*}})}{\theta
<N_{1^{*}}> m(N_{1^{*}})...<N_{4^{*}}> m(N_{4^{*}})} \left(
\frac{N_{3^{*}}}{N_{1^{*}}} \right)^{\frac{1}{2}} \left(
\frac{N_{4^{*}}}{N_{2^{*}}} \right)^{\frac{1}{2}}
N_{5^{*}}^{\frac{3}{2}} N_{6^{*}}^{\frac{3}{2}} \prod_{j=1}^{6}
\| \langle D \rangle I w_{j}(t_{j}) \|_{L_{x}^{2}} \\
& \lesssim \frac{N_{1^{*}}^{--} N_{6^{*}}^{+} }{\theta N^{3-}}
\prod_{j=1}^{6} \| \langle D \rangle I w_{j}(t_{j}) \|_{L_{x}^{2}}
\end{array}
\end{equation}
with $q$ being the logarithmic weight introduced in
(\ref{Eqn:DefWeight}).

\end{document}